\documentclass[graybox]{svmult}

\usepackage[applemac]{inputenc}
\usepackage{amsmath,amssymb}

\numberwithin{equation}{section}

\def\RB{\mathbb{R}}

\def\FB{\mathbb{F}}

\def\FC{\mathcal{F}}

\def\AC{\mathcal{A}}

\def\RC{\mathcal R}

\def\1B{\text{1\!\!I}}

\def\tN{\tilde{N}}

\def\hu{\hat{u}}

\def\hy{\hat{y}}
\def\hk{\hat{k}}
\def\hz{\hat{z}}

\begin{document}
\title*{A stochastic HJB equation for optimal control of forward-backward SDEs}

\date{23 December 2014}

\author{
Bernt \O ksendal and Agn\`es Sulem and Tusheng Zhang}

\institute{Bernt \O ksendal \at \O ksendal, Dept. of Mathematics, University of Oslo, and
Norwegian School of Economics, Helleveien
30, N--5045 Bergen, Norway. 
The research leading to these results has received funding from the
European Research Council under the European Community's Seventh
Framework Programme (FP7/2007-2013) / ERC grant agreement no [228087], P.O. Box 1053 Blindern, N--0316 Oslo, Norway,\email{oksendal@math.uio.no}, \and
Agn\`es Sulem \at Sulem, NRIA Paris-Rocquencourt, Domaine de Voluceau, Rocquencourt, BP 105, Le Chesnay Cedex, 78153, France, and Université Paris-Est, F-77455 Marne-la-Vallée, France, and Dept. of Mathematics, University of Oslo,
P.O. Box 1053 Blindern, N--0316 Oslo, Norway \email{agnes.sulem@inria.fr}, \and
Tusheng Zhang \at Zhang, School of Mathematics, University of Manchester, Oxford Road, Manchester M139PL, United Kingdom, \email{Tusheng.zhang@manchester.ac.uk} }

\maketitle

\titlerunning{A stochastic HJB equation for optimal control of forward-backward SDEs}
\authorrunning{Bernt \O ksendal and Agn\`es Sulem and Tusheng Zhang}
\abstract{
We study optimal stochastic control problems of general coupled systems of forward-backward stochastic differential equations with jumps. By means of the Itô-Ventzell formula the system is transformed into a controlled partial backward stochastic differential equation (PBSDE) with jumps. Using a comparison principle for such BSPDEs we obtain a general stochastic Hamilton-Jacobi- Bellman (HJB) equation for such control problems. In the classical Markovian case with optimal control of jump diffusions, the equation reduces to the classical HJB equation.
The results are applied to study risk minimization in financial markets.
}
\newpage

\section{Introduction}\label{intro}
This paper deals with an HJB equation approach to optimal control of coupled systems of non-Markovian forward-backward stochastic differential equations with jumps. The background for our paper can be divided into two parts:
\begin{itemize}
\item
 \emph{Optimal control of stochastic differential equations (SDEs).}\\

In classical theory of stochastic control of systems described by a stochastic differential equations  there are two important solution methods:\\

(a) \emph{Dynamic programming}, which leads to the classical Hamilton-Jacobi-Bellman (HJB) equation. This is a deterministic non-linear partial differential equation (PDE) in the (unknown) value function for the problem.\\

(b) \emph{The maximum principle}, which involves the maximization of the Hamiltonian and an associated backward stochastic differential equation (BSDE) in the (unknown) adjoint processes. \\

Dynamic programming is a very efficient solution method, but it only works if the system is Markovian. The maximum principle, on the other hand, works also in non-Markovian settings, but the drawback is that it leads to a coupled system of forward-backward SDEs (FBSDEs) with constraints, and this system is difficult to solve in general.
In view of this it is natural to ask if there is an extension of the HJB approach to non-Markovian systems.
The answer has been known to be yes for some time, at least in some cases. See e.g. \cite{BM} and \cite{P}. In \cite{P} a stochastic version of the classical HJB equation is obtained, in the form of a \emph{partial} backward stochastic differential equation (PBSDE), and existence and uniqueness is proved for this type of PBSDEs. However, there it is assumed that the control does not enter the diffusion coefficient of the controlled SDE, and it is assumed that the SDE is driven by Brownian motion only. \\

\item
 \emph{Coupled systems of forward-backward stochastic differential equations (FBSDEs).}\\

It is well-known that in many cases the solution of a coupled system of FBSDEs can be expressed in terms of a solution of a \emph{partial} BSDE.  See e.g.  \cite{MPY} for the Markovian case (which leads to a deterministic backward PDE). For the more general, possibly non-Markovian case, which leads to a PBSDE, see e.g. \cite{MYZ} and the references therein. None of these papers deal with control.
\end{itemize}

The purpose of this paper is to extend the result of \cite{P} above in several directions, in the sense that we obtain a stochastic HJB equation for optimal control in the following context:\\
(i) We consider optimal control of a non-Markovian \emph{coupled system of FBSDEs}.\\
(ii) The system is driven by both a Brownian motion and a Poisson random measure. \\
(iii) The control is allowed to appear in all the coefficients of the system.\\

Our method is based on an extension of the connection in \cite{MYZ}, combined with comparison principles for PBSDEs.
 If the system is a Markovian SDE, then our stochastic HJB equation becomes deterministic and coincides with the classical HJB equation.\\
 
In the last part of the paper we illustrate our theory by studying some applications to finance. In particular,  we apply our results to study a problem of risk minimization in a non-Markovian financial market. \\

%The novelty of our paper lies in the application of the results in \cite{MYZ} to optimal control of FBSDEs and in the extension to jump models.

%SECTION 2
\section{Optimal control of FBSDEs}\label{sec2}
We refer to \cite{OS1} for information about stochastic calculus and control for jump diffusions.
Consider the following controlled coupled FBSDE:
The forward equation in $X(t)$ has the form\\
\begin{equation}\label{eq1.1}
\begin{cases}
dX(t) &= \alpha(t,X(t),Y(t),Z(t),K(t,\cdot), u(t,X(t)))dt \\
 & + \beta (t,X(t),Y(t), Z(t), K(t,\cdot), u(t,X(t)))dB(t) \\
 & + \displaystyle \int_{\RB_0} \gamma(t,X(t),Y(t),Z(t),K(t,\cdot),u(t,X(t)),\zeta) \tN(dt,d\zeta) \; ; \; t \in [0,T] \\
 X(0) & = x \in \RB
 \end{cases}
 \end{equation}
and the backward equation in $Y(t), Z(t), K(t,\zeta)$ has the form
 \begin{align}\label{eq1.2}
% \begin{cases}
 dY(t) & = - g(t,X(t), Y(t), Z(t), K(t,\cdot), u(t,X(t)))dt + Z(t) dB(t) \nonumber \\
 & \displaystyle + \int_{\RB_0} K(t,\zeta) \tN(dt,d\zeta) \; ; \; t \in [0,T];  \; 
 Y(T)  = h(X(T)).
% \end{cases}
 \end{align}
 Here $B(t) = B(t,\omega)$ and $\tN(dt, d\zeta) = N(dt,d\zeta) - \nu(d \zeta)dt \; ; \; t \in [0,T], \omega \in \Omega, \zeta \in \RB_0 := \RB - \{0\}$ is a Brownian motion and an (independent) compensated Poisson random measure, respectively, on a given filtered probability space $(\Omega, \FC, \FB:=\{\FC_t\}_{t \geq 0}, P)$. The measure $\nu$ is the L\' evy measure of the Poisson random measure $N(\cdot, \cdot)$.
 The given functions
 \begin{align*}
 \alpha(t,x,y,z,k,u,\omega)& : [0,T] \times \RB \times \RB \times \RB \times \RC \times V \times \Omega \rightarrow \RB \\
 \beta(t,x,y,z,k,u,\omega)& : [0,T] \times \RB \times \RB \times \RB \times \RC \times V \times \Omega \rightarrow \RB \\
 \gamma(t,x,y,z,k,u, \zeta,\omega)& : [0,T] \times \RB \times \RB \times \RB \times \RC \times V  \times \RB_0 \times \Omega \rightarrow \RB \\
 g(t,x,y,z,k,u,\omega)& : [0,T] \times \RB \times \RB \times \RB \times \RC \times V \times \Omega \rightarrow \RB
 \end{align*}
 are assumed to be $\FB$-predictable for each $x,y,z,k,u$. $\RC$ denotes the set of functions $k(\zeta) : \RB_0 \rightarrow \RB$ and $V$ is a given set of admissible control values $u(t,x,\omega)$, where $u(t) = u(t,X(t),\omega)$ is our control process. The function $h(x,\omega) : \RB \times \Omega \rightarrow \RB$ is assumed to be $\FC_T$-measurable for each $x$.
 We let $\AC$ denote a given set of controls, contained in the set $\mathcal{A}_0$ of predictable control processes $u(t)$ such that the system 
 \eqref{eq1.1}- \eqref{eq1.2} has a unique solution. For results about existence and uniqueness of FBSDEs see \cite{HP}. A possible choice of $\mathcal{A}$ is the set of all processes $u \in \mathcal{A}_0$ such that 
 \begin{equation}
E[\int_0^Tu^2(t)dt] < \infty.
 \end{equation}
 
  If $u \in \AC$, we call $u$ \emph{admissible}.
 We want to find $\hu \in \AC$ such that
 \begin{equation}\label{eq1.3}
 \sup_{u \in \AC} Y^u(0) = Y^{\hu}(0).
 \end{equation}
First we try to write the solution $Y(t)$ of \eqref{eq1.2} of the form
 \begin{equation}\label{eq1.4}
Y(t) = y(t,X(t))\end{equation}
for some random field $y(t,x) = y(t,x,\omega)$ which, together with $z(t,x)$ and $k(t,x,\zeta)$, satisfies a PBSDE of the form
\begin{equation}\label{eq1.5}
\begin{cases}
dy(t,x) & = A_u(y(\cdot), z(\cdot), k(\cdot))(t,x)dt + z(t,x) dB(t) \\
 & \displaystyle + \int_\RB k(t,x,\zeta) \tN(dt, d\zeta) \; ; \; t \in [0,T] \\
 y(T,x) & = h(x),
 \end{cases}
 \end{equation}
 for some partial integro-differential operator $A_u$ acting on $x$. 
  \begin{remark}
We interpret the equation \eqref{eq1.5} in the weak (variational) sense, which means that 
  $y\in C([0,T]; L^2(\RB))\cap L^2([0,T]; V)$ and for $\phi\in C_0^{\infty}(D) \; ; \; t \geq 0$,
\begin{eqnarray}\label{2.0}
<y(t,\cdot), \phi> &=& <y_0(\cdot),\phi> +\int_0^t<A_u y(s, \cdot), \phi>ds 
%+ \int_0^t<b(s,\cdot,Y(s,\cdot)),\phi>_Kds
 \nonumber\\
 %\quad
&&+\int_0^t <z(s,\cdot), \phi> dB(s) \nonumber\\
&&+ \int_0^t \int_\RB <k(s,\cdot,\zeta),\phi> \tN (ds, d\zeta) ,
\end{eqnarray}
where %$A^*$ is the adjoint operator of $A$, and 
$< , >$ denotes the dual pairing between the space $V$ and its dual $V^*$, where $V:=W_0^{1,2}(D)$ is the Sobolev space of order one with zero boundary condition. Note that with this framework the Itô calculus can be applied to \eqref{eq1.4}. See \cite{P},\cite{PR}.
\end{remark}
 By the It\^o-Ventzell formula (see \cite{OZ} and the references therein), 
 \begin{align}\label{eq1.6}
 dY(t) & = A_u(y(\cdot), z(\cdot), k(\cdot))(t,X(t))dt  %\nonumber \\
 % & \quad 
  + z(t, X(t))dB(t) + \int_\RB k(t,X(t), \zeta) \tN(dt,d\zeta) \nonumber \\
  & + y'(t,X(t)) [\alpha(t)dt + \beta(t)dB(t)] + \frac{1}{2} y''(t,X(t))\beta^2(t)dt \nonumber \\
  & + \int_\RB \{ y(t,X(t) + \gamma(t,\zeta)) - y(t,X(t)) - y'(t,X(t)) \gamma(t,\zeta) \} \nu(d \zeta)dt \nonumber \\
  & + \int_\RB \{ y(t,X(t) + \gamma(t,\zeta)) - y(t,X(t))\} \tN(dt,d\zeta)     + z'(t,X(t)) \beta(t)dt \nonumber \\
    & + \int_\RB \{ k(t,X(t) + \gamma(t,\zeta),\zeta) - k(t,X(t),\zeta)\} \nu(d\zeta)dt \nonumber \\
    & + \int_\RB k(t, X(t^-) + \gamma(t,\zeta),\zeta) \tN(dt,d\zeta),
    \end{align}
where $\displaystyle y'(t,x) = \frac{\partial y}{\partial x}(t,x)$ etc.
and where we have used the shorthand notation
    $\alpha(t) = \alpha(t,X(t),Y(t),Z(t),K(t,\cdot),u(t)) \text{ etc.}$
    Rearranging the terms we see that
\begin{align}\label{eq1.7}
 dY(t) & = [A_u(y(\cdot), z(\cdot), k(\cdot))(t,X(t)) + y'(t,X(t)) \alpha(t) + \frac{1}{2} y''(t,X(t))\beta^2(t)\nonumber \\
  & + \int_\RB \{ y(t,X(t) + \gamma(t,\zeta)) - y(t,X(t)) - y'(t,X(t)) \gamma(t,\zeta) \} \nu(d \zeta)
   + z'(t,X(t)) \beta(t)\nonumber \\
    & + \int_\RB \{ k(t,X(t) + \gamma(t,\zeta),\zeta) - k(t,X(t),\zeta)\} \nu(d\zeta)]dt \nonumber \\
    & + [z(t,X(t)) + y'(t,X(t))\beta(t)]dB(t) \nonumber \\
    & + \int_\RB \{y(t,X(t)+ \gamma(t,\zeta)) - y(t,X(t)) 
     + k(t,X(t) + \gamma(t,\zeta),\zeta)\} \tN(dt,d\zeta).
    \end{align}
Comparing \eqref{eq1.7} with \eqref{eq1.2} we deduce the following theorem
 \begin{theorem}\label{th2.1}
 Suppose that $(y(t,x), z(t,x), k(t,x,\cdot))$  satisfies the PBSDE
 \begin{equation}\label{eq1.11}
 %\begin{cases}
 dy(t,x)  = - A_u(t,x)dt + z(t,x) dB(t) + \int_\RB k(t,x,\zeta) \tN(dt,d\zeta); \;  \; 
 y(T,x)  = h(x)
% \end{cases}
 \end{equation}
 where
     \begin{align}\label{eq1.12}
 A_u(t,x) &= A_u(y(\cdot), z(\cdot), k(\cdot))(t,x)  \nonumber \\
 & :=  g(t,x,y(t,x), z(t,x) + y'(t,x) \beta(t), y(t,x + \gamma(t,\cdot)) - y(t,x)\nonumber\\
 & + k(t,x + \gamma(t,\cdot),\cdot),
 u(t,x))\nonumber \\
  & + y'(t,x)\alpha(t) + \frac{1}{2} y''(t,x) \beta^2(t) + z'(t,x) \beta(t) \nonumber \\
  & + \int_\RB \{y(t,x + \gamma(t,\zeta)) - y(t,x) - y'(t,x) \gamma(t,\zeta)\} \nu(d\zeta) \nonumber \\
  & \quad + \int_\RB \{k(t,x + \gamma(t,\zeta),\zeta) - k(t,x,\zeta)\} \nu(d\zeta).
  \end{align}
  Then $(Y(t), Z(t), K(t,\zeta))$, given by
  \begin{equation}
  Y(t) := y(t,X(t)),
  \end{equation}
  \begin{equation}  
  Z(t) := z(t,X(t)) + y'(t,X(t)) \beta(t),
  \end{equation}
  \begin{equation}
  K(t,\zeta) := y(t,X(t) + \gamma(t,\zeta)) - y(t,X(t)) + k(t,X(t)+\gamma(t,\zeta),\zeta),
  \end{equation}
  is a solution of the FBSDE system \eqref{eq1.1}-\eqref{eq1.2}.
  \end{theorem}

%\begin{align}\label{eq1.8}
%A&(y(\cdot), z(\cdot), k(\cdot))(t,X(t)) \nonumber \\
% & = - [g(t,X(t), Y(t), Z(t), K(t,\cdot),u(t,X(t))) + y'(t,X(t)) \alpha(t) + \frac{1}{2} y''(t,X(t)) \beta^2(t) \nonumber \\
% & + \int_\RB \{y(t,X(t) + \gamma(t,\zeta)) - y(t,X(t)) - y'(t,X(t)) \gamma(t,\zeta)\} \nu (d\zeta) \nonumber \\
% & + z'(t,X(t)) \beta(t) + \int_\RB \{ k(t,X(t) + \gamma(t,\zeta),\zeta) - k(t,X(t),\zeta)\} \nu (d\zeta)]
% \end{align}
% \begin{equation}\label{eq1.9}
% Z(t) = z(t,X(t)) + y'(t,X(t)) \beta(t)
% \end{equation}
% \begin{equation}\label{eq1.10}
% K(t,\zeta) = y(t,X(t) + \gamma(t,\zeta)) - y(t,X(t)) + k(t,X(t)+\gamma(t,\zeta),\zeta).
% \end{equation}
%n the following we will use the shorthand notations
%\begin{equation}
%\beta(t) := \beta(t,x,y(t,x),z(t,x),k(t,x,\cdot)) \text{  etc},
%\end{equation}
%\begin{equation}
%\tz(t,x) := z(t,x) + y'(t,x) \beta(t)
%\end{equation}
%\begin{equation}
% \tk(t,x,\cdot) := y(t,x + \gamma(t,\cdot)) - y(t,x) + k(t,x + \gamma(t,\cdot),\cdot).
% \end{equation}

  \begin{definition}
We say that the PBSDE \eqref{eq1.11} satisfies the comparison principle with respect to $u$ if  for all $u_1, u_2\in \AC$ and all ${\cal F}_T$-measurable random variables $ h_1(x), h_2(x)$ with corresponding solutions $(y_i, z_i, k_i), i=1,2$, of (\ref{eq1.11}) such that 
$$A_{u_1}(t,x)\leq A_{u_2}(t,x)  \text{ for all }  t, x\in [0,T]\times \RB$$
 $$ \text{and } \quad h_1(x)\leq h_2(x)   \text{ for all } x\in \RB,$$
 we have 
 $$y_1(t,x)\leq y_2(t,x) \text{ for all } t, x\in [0,T]\times \RB.$$
\end{definition}
Sufficient conditions for the validity of comparison principles for PBSDEs with jumps is still an open question in this setting. For related results see \cite{OSZ2}. 
However in the Brownian case, sufficient conditions for the validity of comparison principles for PBSDEs of the type 
\eqref{eq1.12} 
are given in Theorem 2.13 in \cite{MYZ}. 
Using this result we get
\begin{theorem}\label{th2.3}
Assume that the following holds:
\begin{itemize}
\item
$N=K=0$, i.e. there are no jumps
\item
The coefficients $\alpha$, $\beta$, and $g$ are $\mathbb{F}$ - progressively measurable for each fixed $(x,y,z)$ and $h(x)$ is $\mathcal{F}_T$ - measurable for each fixed $x$
\item
$\alpha, \beta, g, h$ are uniformly Lipschitz-continuous in $(x,y,z)$
\item
$\alpha$ and $\beta$ are bounded and
$%\begin{equation}
E[\int_0^T g^2(t,0,0,0)dt + h^2(0)] < \infty
$%\end{equation}
\item
$\alpha(t,x,y,z,u)$ does not depend on $z$.
\end{itemize}
Then  the comparison principle holds for the PBSDE \eqref{eq1.11}.
\end{theorem}
  From the above we deduce the following result, which may be regarded as a \emph{stochastic HJB equation} for optimal control of possibly non-Markovian FBSDEs. 
  \begin{theorem}\label{th1.1} (Stochastic HJB equation.)
  Suppose  the comparison principle holds for the BSPDE \eqref{eq1.11}. 
  Moreover, suppose that for all $t,x,\omega$ there exists a maximizer %\linebreak 
  $u = \hu(t,x) = \hu(y,y',y'',z,z',k)(t,x,\omega)$ of the function $u \rightarrow A_u(t,x)$.
 Suppose the system \eqref{eq1.11} with $u = \hu$ has a unique solution $(\hy(t,x),\hz(t,x),\hk(t,x,\cdot))$ and that $\hu(t,X(t)) \in \AC$. Then $\hu(t,X(t))$ is an optimal control for the problem \eqref{eq1.3}, with optimal value
  \begin{equation}\label{eq1.13}
  \sup_{u \in \cal{A}}Y^{u}(0) =Y^{\hu}(0) = \hy(0,x).
  \end{equation}
\end{theorem}
Note that in this general non-Markovian setting the classical \emph{value function} from the dynamic programming is replaced by the solution $\hy(t,x)$ of the PBSDE \eqref{eq1.11} for $u = \hu$. 
%See Example \ref{exa2.1} and Example \ref{exa5.6} below.
\section{Applications}\label{sec4}
We now illustrate Theorem \ref{th1.1} by looking at some examples. First we consider the classical Merton problem. The solution of this problem is well known in the Markovian case with deterministic coefficients, but we consider here  the general non-Markovian case, when the coefficients are stochastic processes:
\begin{example} [Maximizing expected utility from terminal wealth]\label{exa2.1}\rm
Consider a financial market consisting of 
a risk free investment, with unit price
%\begin{equation} \label{eq3.1}
$S_0(t) := 1\; ; \; t\in [0,T],$
%\end{equation}
and a risky investment, with unit price
\begin{equation} \label{eq3.2}
dS_1(t) = S_1(t) [b(t)dt + \sigma(t) dB(t)] \; ; \; t\in [0,T]. 
\end{equation}
Here $b(t) = b(t,\omega)$ and $\sigma(t) = \sigma(t,\omega) > 0$ are given adapted processes.
Let $u(t,X(t))$ be a \emph{portfolio}, representing the \emph{amount} invested in the risky asset at time t. If $u$ is self-financing, then the corresponding wealth $X(t)$ at time $t$ is given by the stochastic differential equation
\begin{equation}\label{eq3.3}
%\begin{cases}
dX(t) = dX_x^{u}(t)=u(t,X(t))[b(t) dt + \sigma(t) dB(t)] \;,  t \in [0,T] \; ; \;
 X(0) = x > 0.
% \end{cases}
 \end{equation}
 Let $(Y(t),Z(t))=(Y_x^{u}(t),Z_x^{u}(t))$ be the solution of the BSDE
 \begin{equation}\label{eq2.2}
% \begin{cases}
 dY(t)  = Z(t) dB(t) \;, \; t \in [0,T] \; ; \; 
  Y(T)  = U(X(T)), 
  %\end{cases}
  \end{equation}
  where $U(X) = U(X,\omega)$ is a given utility function, possibly random. Then
  $$Y_x^{u}(0) = E[U(X_x^{u}(T))].$$
 Therefore, the classical portfolio optimization problem of Merton is to find $\hu \in \AC$ such that
 \begin{equation}\label{eq2.2a}
 \sup_{u \in \AC} Y_x^u(0) = Y_x^{\hu}(0).
 \end{equation}
 In the following we assume that 
 \begin{equation} \label{f}
 \sup_{u \in \AC} Y_x^u(0) < \infty. \end{equation} 
In this general non-Markovian setting with stochastic coefficients $b(t)=b(t,\omega)$ and $\sigma(t)=\sigma(t,\omega)>0$, an explicit expression for the optimal portfolio $\hu$ is not known.
  We  apply the theory from the previous sections to study this problem. 
  In this case we get, from \eqref{eq1.12},
  \begin{equation}\label{eq2.3}
  A_u(t,x) = y'(t,x) u b(t) + \frac{1}{2} y''(t,x) u^2 \sigma^2(t,x) + z'(t,x) u \sigma(t)
  \end{equation}
  which is maximal when
  \begin{equation}\label{eq2.4}
  u = \hu(t,x) = - \frac{y'(t,x) b(t) + z'(t,x) \sigma(t)}{y''(t,x) \sigma^2(t)}.
  \end{equation}
  Substituting this into $A_{\hu}(t,x)$ we obtain
  \begin{equation}\label{eq2.5}
  A_{\hu}(t,x) = - \frac{(y'(t,x)b(t) + z'(t,x)\sigma(t))^2}{2y''(t,x) \sigma^2(t)}.
  \end{equation}
%\end{example}
Hence the PBSDE for $y(t,x)$ gets the form
\begin{equation}\label{eq2.6}
\begin{cases}
dy(t,x) & = \displaystyle \frac{(y'(t,x)b(t) + z'(t,x)\sigma(t))^2}{2y''(t,x)\sigma^2(t)} dt + z(t,x) dB(t) \; ; \; t \in [0,T] \\
y(T,x) & = U(x).
\end{cases}
\end{equation}
We have proved:
\begin{proposition} \label{th3.1}
Suppose there exists a solution $(y(t,x), z(t,x))$ of the PBSDE \eqref{eq2.6} with $y''(t,x)<0$. Suppose that $\hu$ defined in \eqref{eq2.4} is admissible. Then $\hu$ is optimal for problem \eqref{eq2.2a} and
\begin{equation} \label{eq2.6a}
y(0,x) =  \sup_{u \in \AC} Y_x^u(0) = Y_x^{\hu}(0).
\end{equation}
\end{proposition}
Note that if $b,\sigma$ and $U$ are {\it deterministic}, we can choose $z(t,x) = 0$ in \eqref{eq2.6} and this leads to the following (deterministic) PDE for $y(t,x)$:
\begin{equation}\label{eq2.7}
%\begin{cases}
\frac{\partial y}{\partial t}(t,x) - \frac{y'(t,x)^2b^2(t)}{2y''(t,x) \sigma^2(t)} = 0 \; ; \; t \in [0,T] \; ; \; 
y(T,x)=U(x).
%\end{cases}
\end{equation}
This is the classical Merton PDE for the value function, usually obtained by dynamic programming and the HJB equation.
Hence we may regard \eqref{eq2.4}-\eqref{eq2.6} as a generalization of the Merton equation \eqref{eq2.7} to the non-Markovian case with stochastic $b(t), \sigma(t)$ and $U(x)$.
%, in which case classical dynamic programming cannot be used.
%Thus we see that
 The Markovian case corresponds to the special case when $z(t,x) = 0$ in the BSDE \eqref{eq2.7}. Therefore $\hy(s,x)$ is a stochastic generalization of the value function
\begin{equation}\label{eq5.23}
\varphi(s,x) := \sup_{u\in \cal{A}} E[U(X^{u}_{s,x}(T))]
\end{equation}
where
\begin{equation}\label{eq5.24}
%\begin{cases}
dX_{s,x}^{u}(t)  = u(t) [b(t) dt + \sigma(t) dB(t)] \; ; \; t \geq s \; ; \; 
X_{s,x}^{u}(s)  = x.
%\end{cases}
\end{equation}
Let us compare with the use of the classical HJB:
\begin{equation}\label{eq3.10}
\begin{cases}
\displaystyle \frac{\partial \varphi}{\partial s} (s,x) + \max_v \left\{ \frac{1}{2} v^2 \sigma_0^2(s) \varphi''(s,x) + v b_0(s) \varphi'(s,x) \right\} & = 0 \; ; \; s < T \\
\varphi(T,x) & = U(x).
\end{cases}
\end{equation}
The maximum is attained at
\begin{equation}\label{eq5.26}
v = \hu(s,x) = - \frac{b_0(s) \varphi'(s,x)}{\varphi''(s,x) \sigma^2_0(s)} .
\end{equation}
Substituted into \eqref{eq3.10} this gives the HJB equation
\begin{equation}\label{eq5.27}
\frac{\partial \varphi}{\partial s}(s,x) - \frac{\varphi'(s,x)^2 b_0^2(s)}{\varphi'' (s,x) \sigma^2_0(s)} = 0,
\end{equation}
which is identical to \eqref{eq2.7}.
\end{example}
\begin{example}[Risk minimizing portfolios]\label{exa5.6}
\rm
Now suppose $X(t) = X^{u}_x(t)$ is as in \eqref{eq3.3}, while $(Y(t),Z(t))=(Y_x^{u}(t),Z_x^{u}(t))$ is given by the BSDE
\begin{equation}\label{eq5.28}
%\begin{cases}
dY(t) =  - ( - \frac{1}{2} Z^2(t))dt + Z(t) dB(t) \; ; \; 
 Y(T) = X(T).
% \end{cases}
 \end{equation}
 Note that the driver $\displaystyle g(z): = - \frac{1}{2} z^2$ is concave.
 We want to minimize the \emph{risk}  of the terminal financial standing $X(T)$, denoted by $\rho(X(T))$. If we interpret the risk in the sense of the \emph{convex risk measure} defined in terms of the BSDE \eqref{eq5.28} we have 
 $$ \rho(X(T)) = -Y(0).$$
 % where $\rho(X(T))$ is the \emph{risk} of $X(T)$ with respect to the driver $\displaystyle g(z) = \frac{1}{2} z^2$. 
 See e.g. \cite{QS}, \cite{R} for more information about the representation of risk measures via BSDEs. 
 Thus the risk minimization problem is to find $\hu \in \AC$ such that
 \begin{equation}\label{eq5.29a}
 \inf_{u \in \AC} -Y_x^u(0) = -Y_x^{\hu}(0),
 \end{equation}
where $Y_x^{u}(t)$ is given by \eqref{eq5.28}.
 By changing sign we can consider the supremum problem in stead. In this case we get
 \begin{align}\label{eq5.29}
 A_u(t,x) & = - \frac{1}{2} (z(t,x) + y'(t,x) u \sigma(t))^2  + y'(t,x) u b(t)\nonumber \\
 &  + \frac{1}{2} y''(t,x) u^2 \sigma^2 (t) + z'(t,x) u \sigma(t),
 \end{align}
 which is minimal when $u = \hu(t,x)$ satisfies
 \begin{equation}\label{eq5.30}
 \hu(t,x) = - \frac{z(t,x) y'(t,x) \sigma(t) - y'(t,x) b(t)  - z'(t,x) \sigma(t)}{((y'(t,x))^2 - y''(t,x)) \sigma^2(t)}.
 \end{equation}
 This gives
 \begin{equation}\label{eq5.31}
 A_{\hu}(t,x) = - \frac{1}{2} \hz^2(t,x) + \frac{(\hz(t,x) \hy'(t,x) \sigma(t) - \hy'(t,x) b(t)  - \hz'(t,x) \sigma(t))^2}{2((\hy'(t,x))^2 - \hy''(t,x)) \sigma^2(t)}.
 \end{equation}
 and hence $(\hy(t,x), \hz(t,x))$ solves the PBSDE
 \begin{equation}\label{eq5.32}
% \begin{cases}
 d\hy(t,x) = - A_{\hu}(t,x) dt + \hz(t,x) dB(t), \; 0 \leq t \leq T \; ; \; 
 \hy(T,x) = x.
% \end{cases}
 \end{equation}
 We have proved: 
 \begin{proposition}
 Suppose there exists a solution $(\hy(t,x),\hz(t,x))$ of the PBSDE \eqref{eq5.32}. Suppose $\hu$ defined by \eqref{eq5.30} belongs to $\AC$. Then $\hu$ is optimal for the risk minimizing problem \eqref{eq5.29a}, and the 
 minimal risk is 
 \begin{equation}
  \inf_{u \in \AC} -Y_x^u(0) = -Y_x^{\hu}(0) = - \hy(0,x).
 \end{equation}
 \end{proposition}
Next we look at the special case when $b(t)$ and $\sigma(t)$ are \emph{deterministic}. 
 Let us try to choose $\hz(t,x) = 0$ in \eqref{eq5.32}.
 Then this PBSDE reduces to the (backward) PDE
 \begin{equation}\label{eq5.33}
 \begin{cases}
 \displaystyle \frac{\partial \hy(t,x)}{\partial t} = - \frac{(\hy'(t,x) b(t))^2}{2((\hy'(t,x))^2 - \hy''(t,x)) \sigma^2(t)} \; ; \; 0 \leq t \leq T \\
 \hy(T,x) = x.
 \end{cases}
 \end{equation}
 We try a solution of the form
 \begin{equation}\label{eq5.34}
 \hy(t,x) = x + a(t),
 \end{equation}
 where $a(t)$ is deterministic. Substituted into \eqref{eq5.33} this gives
 \begin{equation}\label{eq5.35}
% \begin{cases}
 a'(t) = \displaystyle - \frac{1}{2} \left( \frac{b(t)}{\sigma(t)}\right)^2 \;, \; 0 \leq t \leq T \; ; \; 
 a(T) = 0
% \end{cases}
 \end{equation}
 which gives
 $$a(t) = \int_t^T \frac{1}{2} \left( \frac{b(s)}{\sigma(s)}\right)^2 ds \; ; \; 0 \leq t \leq T.$$
 With this choice of $a(t)$,   \eqref{eq5.33} is satisfied and we  conclude that the minimal risk is
 \begin{equation}\label{eq5.36}
 \rho_{min}(X(T)) = - Y^{(\hu)}(0) = - \hy(0,x) = - x - \int_0^T \frac{1}{2} \left( \frac{b(s)}{\sigma(s)}\right)^2 ds
 \end{equation}
 Hence by \eqref{eq5.30} the optimal (risk minimizing) portfolio is
 \begin{equation}\label{eq5.37}
 \hu(t,X(t)) = \frac{ b(t)}{\sigma^2(t)}.
 \end{equation}
\begin{remark} Note that \eqref{eq5.36} can be interpreted by means of \emph{entropy} as follows:
 Recall that in general the {\it entropy} of a measure $Q$ with respect to the measure $P$ is defined by 
 $$\displaystyle H(Q \mid P) := E \left[ \frac{dQ}{dP} \ln \frac{dQ}{dP}\right].$$
 Define
 \begin{equation} 
 \Gamma(t) = \exp \left(- \int_0^t \frac{b(s)}{\sigma(s)} dB(s) - \frac{1}{2} \int_0^t ( \frac{b(s)}{\sigma(s)})^2 ds \right).
 \end{equation}
% Then
% \begin{equation}
% d \Gamma(t) = \Gamma(t) \left[ -\frac{b(t)}{\sigma(t)} dB(t)\right]  \; ; \; \Gamma(0) = 1.
% \end{equation}
 By the It\^o formula we have
 \begin{align*}
 d(\Gamma(t) \ln \Gamma(t)) &= \Gamma(t) \left[ - \frac{b(t)}{\sigma(t)} dB(t) - \frac{1}{2} \left( \frac{b(t)}{\sigma(t)} \right)^2 dt\right] \\
 & + (\ln \Gamma(t)) \Gamma(t) \left( - \frac{b(t)}{\sigma(t)} dB(t)\right) + \Gamma(t) \left( - \frac{b(t)}{\sigma(t)} \right) \left( -\frac{b(t)}{\sigma(t)} \right)dt.
 \end{align*}
 Hence, if we define the measure $Q_{\Gamma}(\omega)$ by
 \begin{equation}
 dQ_{\Gamma}(\omega) := \Gamma(T) dP(\omega)
 \end{equation}
 we get
 \begin{align*}
 E &\left[ \frac{dQ_{\Gamma}}{dP} \ln \frac{dQ_{\Gamma}}{dP}\right] = E [ \Gamma(T) \ln \Gamma(T)] \\
 & = E \left[ \int_0^T \Gamma(t) \frac{1}{2} \left( \frac{b(t)}{\sigma(t)}\right)^2 dt \right] =  \frac{1}{2}\int_0^T \left( \frac{b(t)}{\sigma(t)}\right)^2 dt,
 \end{align*}
 which proves that \eqref{eq5.36} can be written
 \begin{equation}
 \rho_{min}(X(T)) = -x - H(Q_{\Gamma} \mid P).
\end{equation}
Note that $Q_{\Gamma}$ is the \emph{unique equivalent martingale measure} for the market % \eqref{eq3.1},
\eqref{eq3.2}. 
\end{remark}
Thus we have proved that if 
%\begin{theorem}
%Suppose that
 the coefficients $b(t)$ and $\sigma(t)$ in \eqref{eq3.3} are deterministic and if the portfolio
$\hu(t,X(t)):= \frac{ b(t)}{\sigma^2(t)}$
is admissible, then $\hu$ is a risk minimizing portfolio for the problem \eqref{eq5.29a} and the minimal risk is
%\begin{equation}
%\inf_{u \in \AC} \rho(X(T)) = - Y_x^{\hu}(0) = - \hy(0,x) = - x - \int_0^T \frac{1}{2} \left( \frac{b(s)}{\sigma(s)}\right)^2 ds
% = -x - H(Q_{\Gamma} \mid P).
%\end{equation}
%Thus we conclude that the minimal risk 
is equal to minus the initial wealth x minus the entropy of the equivalent martingale measure. 
%for the market \eqref{eq3.2}.
%\end{theorem}
%It is natural to ask what the corresponding result would be in the incomplete market case.
For alternative solution approaches to this problem based on 
(i) the maximum principle for optimal control of forward-backward SDEs, and on 
(ii) stochastic differential games,
see the survey paper \cite{OS2}.
\end{example}

\end{document}